\documentclass{amsart}
\usepackage{amsmath,amssymb,amsbsy,amsfonts,amsthm,latexsym,
                       amsopn,amstext,amsxtra,euscript,amscd}
                       \usepackage{color}
\usepackage{array}
\usepackage{ifthen}
\usepackage{url}

\newtheorem{theorem}{Theorem}[section]
\newtheorem{lemma}{Lemma}[section]

\theoremstyle{definition}

\def\L{\mathbb L}
\def\R{\mathbb R}

\def\hh{{\rm h}\kern.4pt}

\begin{document}

\title[On terms of Pell equation with Tribonacci numbers]
{On the $X$-coordinates of Pell equations which are Tribonacci numbers}

\author[Florian Luca]{Florian Luca}
\address{School of Mathematics,  University of the Witwatersrand,  Private Bag X3, Wits 2050, South Africa and Centro de Ciencias Matem\'aticas UNAM, Morelia, Mexico}
\email{florian.luca@wits.ac.za}

\author[Amanda Montejano]{Amanda Montejano}
\address{Facultad de Ciencias, UNAM Campus Juriquilla, Mexico}
\email{amandamontejano@ciencias.unam.mx}

\author[Laszlo Szalay]{Laszlo Szalay}
\address{Department of Mathematics and Informatics, J.~Selye University, Hradna ul.~21, 94501 Komarno, Slovakia and
Institute of Mathematics, University of West Hungary, 9400, Sopron, Ady \'ut 5, Hungary}
\email{szalay.laszlo@nyme.hu}

\author[Alain Togb\'e]{Alain Togb\'e}
\address{Department of Mathematics, Statistics and Computer Science, Purdue University Northwest, 1401 S, U.S. 421, Westville IN 46391 USA}
\email{atogbe@pnw.edu}

\date{\today}


\begin{abstract}
For an integer $d\geq 2$ which is not a square, we show that there is at most one value of the positive integer $X$ participating in the Pell equation $X^2-dY^2=\pm 1$ which is a Tribonacci number, with a few exceptions that we completely characterize.
\end{abstract}

\maketitle

\section{Introduction}\label{sec:1}

Let $d>1$ be a positive integer which is not a perfect square. Consider the Pell equation
\begin{equation} \label{eq1}
 X^2-dY^2=\pm 1.
\end{equation}
All its positive integer solutions $(X,Y)$ are given by
$$ 
X_n+ Y_n {\sqrt{d}} =(X_1+Y_1{\sqrt{d}})^n
$$
for some positive integer $n$, where $(X_1,Y_1)$ is the smallest positive solution. In several recent papers, the following problem was investigated. Let ${\bf U}=\{U_n\}_{n\ge 0}$ be some interesting sequence of positive integers. 
What can one say about the square-free integers $d$ such that the equation $X_n\in {\bf U}$ has at least two solutions $n$? For most sequences, one expects that the answer to such a question would be 
that the equation $X_n\in {\bf U}$ has at most one  positive integer solution $n$ for any given $d$ except maybe for a few (finitely many) values of $d$. In \cite{DLT}, this was shown to be so when ${\bf U}$ is the sequence of all 
base $10$-repdigits; that is, numbers of the form $c(10^m-1)/9$, for some positive integers $m\ge 1$ and $c\in \{1,\ldots,9\}$.  The only exceptional $d$'s in this case were $d=2,3$. For each of these two values of $d$, the equation $X_n\in {\bf U}$ 
has two solutions $n$. In \cite{FL}, it was shown, more generally, that if $b\ge 2$ is any fixed positive integer, ${\bf U}$ is the sequence of base $b$-repdigits,  and $d$ is such that $X_n\in {\bf U}$ has two solutions $n$, then 
$$
d<\exp((10b)^{10^5}).
$$  
In \cite{LT}, it was shown that if ${\bf U}$ is the sequence of Fibonacci numbers, then the equation $X_n\in {\bf U}$ has at most one positive integer solution $n$, except when $d=2$ for which there are exactly two solutions.

In this paper, we consider the same problem for the sequence ${\bf U}:={\bf T}$ of Tribonacci numbers given by $T_0=0,~T_1=T_2=1$ and 
$T_{m+3}=T_{m+2}+T_{m+1}+T_m$, for all $m\ge 0$. Our result is the following.

\begin{theorem}\label{main}
Let $d\ge 2$ be square-free. The Diophantine equation
 \begin{equation}\label{eq4}
X_n=T_m
\end{equation}
has at most one solution $(n, m)$ in positive integers with the following exceptions:

\begin{itemize}
	\item $(n_1,m_1)=(1,3)$ and $(n_2,m_2)=(2,5)$ in the $+1$ case,
	\item $(n_1,m_1)=(1,1)$, $(n_2,m_2)=(1,2)$ and $(n_3,m_3)=(3,5)$ in the $-1$ case. 	
\end{itemize}

\end{theorem}

A few words about our method. For the arguments in \cite{DLT}, \cite{FL} and \cite{LT}, the arithmetical properties of the members of ${\bf U}$ played an important role. For example, it was important to know all the solutions of the equation $U_m=2X^2-1$ 
in positive integers $(m,X)$, which are easy to find when ${\bf U}$ is the sequence of Fibonacci numbers or base $10$-repdigits. It was also important that $\gcd(U_m,U_n)$ was closely related to $U_{\gcd(m,n)}$, which is the case both when ${\bf U}$ is the sequence of Fibonacci numbers and the sequence of repdigits. In contrast, the sequence of Tribonacci numbers does not display similar properties. For example, the equation $T_m=2X^2-1$ in positive integers $(m,X)$ is unsolved and there is no general method that would allow one to solve such equation (albeit some tricky elementary arguments might solve such equation) and $\gcd(T_m,T_n)$ is not related in any obvious way to $T_{\gcd(m,n)}$. Our method consists in applying Baker's 
theory of linear forms in logarithms three times to three different linear forms in order to get an absolute bound on all the variables, after which we use reduction procedures to reduce our bounds to some reasonable values and carry on the computations in the remaining range.   
Our method works equally well not only for the Tribonacci sequence but for other linearly recurrent sequences satisfying certain technical conditions.   For example, it works for sequences $(u_m)_{m\ge 1}$ which are linearly recurrent, nondegenerate, have a simple dominant root
$\alpha>1$ and all other roots of absolute value smaller than $1$, and furthermore if $a$ is the coefficient of $\alpha^m$ in the Binet formula for $u_m$, then $\log (2a)$ and $\log \alpha$ are linearly dependent over ${\mathbb Q}$ (which insures that the analog of the left--hand side of \eqref{eq:2logs} is nonzero). 

\section{The Tribonacci sequence}\label{sec:2}

Here, we recall a few important properties of the Tribonacci sequence $\{T_n\}_{n\ge 0}$.
The characteristic equation
$$
x^3-x^2-x-1=0
$$
has roots $\alpha,~\beta,~\gamma={\overline{\beta}}$, where
$$
\alpha=\frac{1+\omega_1+\omega_2}{3},\qquad \beta=\frac{2-\omega_1-\omega_2+{\sqrt{3}} i (\omega_1-\omega_2)}{6},
$$
and
$$
\omega_1=\sqrt[3]{19+3{\sqrt{33}}}\quad {\text{\textrm{and}}}\quad \omega_2=\sqrt[3]{19-3{\sqrt{33}}}.
$$
Further, Binet's formula is
\begin{equation}
\label{eq:Binet}
T_m=a\alpha^m+b\beta^m+c\gamma^m,\quad {\text{\textrm{for all}}}\quad m\ge 0,
\end{equation}
where
\begin{equation}
\label{eq:coefficients}
a=\frac{1}{(\alpha-\beta)(\alpha-\gamma)},\quad b=\frac{1}{(\beta-\alpha)(\beta-\gamma)},\quad c=\frac{1}{(\gamma-\alpha)(\gamma-\beta)}={\overline{b}}
\end{equation}
(see \cite{S}). Numerically,
\begin{equation}
\label{eq:num}
\begin{split}
& 1.83<\alpha<1.84,\\
& 0.73<|\beta|=|\gamma|=\alpha^{-1/2}<0.74,\\
& 0.18<a<0.19,\\
& 0.35<|b|=|c|<0.36.
\end{split}
\end{equation}
Further,
\begin{equation}
\label{eq:BL}
\alpha^{m-2}\le T_m\le \alpha^{m-1},
\end{equation}
for all $m\ge 2$ (see \cite{BL}).

\section{Linear forms in logarithms}\label{sec:3}

We need some results from the theory of lower bounds in nonzero linear forms in logarithms of algebraic numbers. We start by recalling Theorem 9.4 of \cite{BMS}, which is a modified version of a result of Matveev \cite{Matveev}.  Let $\L$ be an algebraic number field of degree $d_{\L}$. Let $\eta_1, \eta_2, \ldots, \eta_l \in \L$  not $0$ or $1$ and $d_1, \ldots, d_l$ be nonzero integers. We put
$$
D =\max\{|d_1|, \ldots, |d_l|, 3\},
$$
and put
$$
\Gamma = \prod_{i=1}^l \eta_i^{d_i} -1.
$$
Let $A_1, \ldots, A_l$ be positive integers such that 
$$
A_j \geq h'(\eta_j) := \max \{d_{\L}h(\eta_j), |\log \eta_j|, 0.16\}, \quad {\text{\rm for}}\quad j=1,\ldots l,
$$
where for an algebraic number $\eta$  of minimal polynomial
$$
f(X)=a_0(X-\eta^{(1)})\cdots(X-\eta^{(k)})\in {\mathbb Z}[X]
$$
over the integers with positive $a_0$, we write $h(\eta)$ for its Weil height given by
$$
h(\eta)=\frac{1}{k}\left(\log a_0+\sum_{j=1}^k \max\{0,\log |\eta^{(j)}|\}\right).
$$
The following consequence of Matveev's theorem is Theorem 9.4 in \cite{BMS}.
\begin{theorem}
\label{thm:Matveev}
 If $\Gamma \neq 0$ and $\L \subseteq \R $, then
\begin{equation*}
\label{ineq:matveev} \log |\Gamma| > -1.4 \cdot
30^{l+3}l^{4.5}d_{\L}^2(1+ \log d_{\L})(1+ \log D)A_1A_2\cdots A_l.
\end{equation*}
\end{theorem}
When $k=2$ and $\eta_1,~\eta_2$ are positive and multiplicatively independent, we can do better. Namely, let in this case $B_1,~B_2$ be real numbers larger than $1$ such that
$$
\log B_i\ge \max\left\{h(\eta_i), \frac{|\log \eta_i|}{d_{\L}}, \frac{1}{d_{\L}}\right\}\qquad i=1,2,
$$
and
$$
b':=\frac{|d_1|}{d_{\L} \log B_2}+\frac{|d_2|}{d_{\L} \log B_1}.
$$
Put
$$
\Lambda=d_1\log \eta_1+d_2\log\eta_2.
$$
Note that $\Lambda\ne 0$ when $\eta_1$ and $\eta_2$ are multiplicatively independent. The following inequality is Corollary 2 in \cite{LMN}.

\begin{theorem}
\label{thm:LMN}
With the above notation, assuming that  $k=2$, $\L$ is real, $\eta_1,~\eta_2$ are positive and multiplicatively independent, then
\begin{equation}
\label{eq:LMN}
\log |\Lambda| >-24.34 d_{\L}^4 \left(\max\left\{\log b'+0.14, \frac{21}{d_{\L}}, \frac{1}{2}\right\}\right)^2\log B_1\log B_2.
\end{equation}
\end{theorem}

\section{The Baker-Davenport lemma}\label{sec:4}

We recall the Baker-Davenport reduction method (see \cite[Lemma 5a]{dujella}), which will be useful to reduce the bounds arising from applying Theorems \ref{thm:Matveev} and \ref{thm:LMN}.

\begin{lemma}\label{lem:Baker-Davenport}
Let $\kappa\ne 0$ and $\mu$ be real numbers. Assume that $M$ is a positive integer. Let $P/Q$ be the convergent
of the continued fraction expansion of $\kappa$ such that $Q > 6M$
and put
$$
\xi=\left\| \mu Q \right\| - M \cdot \| \kappa Q\|,
$$
where $\parallel\cdot\parallel$ denotes the distance from the
nearest integer. If $\xi>0$, then there is no solution of the
inequality
$$0 < |m\kappa -n+\mu| < AB^{-k}$$
in positive integers $m$, $n$ and $k$ with
$$
\frac{\log\left(AQ/\xi\right)}{\log B}\leq k  \qquad and\qquad m\le M.
$$
\end{lemma}

\section{Bounding the variables}\label{sec:5}

We assume that $(X_1,Y_1)$ is the minimal solution of the Pell equation \eqref{eq1}. Setting
$$
X_1^2-dY_1^2=:\varepsilon,\qquad \varepsilon\in \{\pm 1\},
$$
we put 
$$
\delta:=X_1+{\sqrt{d}} Y_1\qquad {\text{\rm and}}\qquad \eta:=X_1-{\sqrt{d}}Y_1=\varepsilon \delta^{-1}.
$$
Then
\begin{equation}
\label{eq:Binet1}
X_n=\frac{1}{2}(\delta^n+\eta^n).
\end{equation}
Since $\delta\ge 1+{\sqrt{2}}$, it follows that the estimate
\begin{equation}
\label{eq:delta}
\frac{\delta^{n}}{\alpha}\le X_n<\delta^n\qquad {\text{\rm holds~for~all}}\qquad n\ge 1.
\end{equation}
We now assume that $(n_1,m_1)$ and $(n_2,m_2)$ are pairs of positive integers such that
$$
X_{n_1}=T_{m_1}\qquad {\text{\rm and}}\qquad X_{n_2}=T_{m_2}.
$$
To fix ideas, we assume that $n_1<n_2$, so $m_1<m_2$. Setting $(n,m):=(n_i,m_i)$, for $i\in \{1,2\}$ and using inequalities \eqref{eq:BL} and \eqref{eq:delta}, we get that
\begin{equation}
\label{eq:ttt}
\alpha^{m-2}\le T_m=X_n< \delta^n\qquad {\text{\rm and}}\qquad \frac{\delta^{n}}{\alpha}\le  X_n=T_m\le \alpha^{m-1}.
\end{equation}
Hence,
\begin{equation}
\label{eq:mn}
nc_1\log \delta\le m\le n c_1\log \delta+2,\qquad c_1:=1/\log \alpha
\end{equation}
holds. Next, using \eqref{eq:Binet} and \eqref{eq:Binet1}, we get
$$
\frac{1}{2} (\delta^n+\eta^n)=a\alpha^m+b\beta^m+c\gamma^m,
$$
so
$$
\delta^n (2a)^{-1} \alpha^{-m}-1=-(2a)^{-1} \alpha^{-m} \eta^n+(b/a) (\beta\alpha^{-1})^m+(c/a)(\gamma \alpha^{-1})^m.
$$
Hence, using \eqref{eq:num}, and assuming that $m>100$, we have
\begin{eqnarray*}
\left| \delta^n (2a)^{-1} \alpha^{-m}-1\right| & \le & \frac{1}{2a \alpha^m \delta^n}+\frac{|b| |\beta|^m}{a \alpha^m}+\frac{|c| |\gamma|^m}{a\alpha^m}\\
&< & \frac{1}{2a \alpha^m \delta^n}+\frac{2|b|}{a \alpha^{3m/2}}\\
& < & \frac{\alpha^3}{2a \alpha^{2m}}+\frac{2|b|}{a \alpha^{3m/2}}\\
& < & \frac{4.5}{\alpha^{3m/2}}.
\end{eqnarray*}
In the above, we used that $|b|/a<2$ (see \eqref{eq:num}) and that $\alpha^{m/2}>\alpha^3/(2a)$ which holds for $m>100$. Since $\alpha^{3m/2}>6$, it follows that the last number above is $<1/2$. Thus,
\begin{equation}
\label{eq:multform}
\left| \delta^n (2a)^{-1} \alpha^{-m}-1\right| <\frac{4.5}{\alpha^{3m/2}}.
\end{equation}
Put
$$
\Lambda:=n\log \delta-\log 2a-m\log \alpha.
$$
Since $|e^{\Lambda}-1|<1/2$, it follows that
$$
|\Lambda|<2|e^{\Lambda}-1|<\frac{9}{\alpha^{3m/2}}.
$$
Recalling that $(m,n)=(m_i,n_i)$, we get that 
\begin{equation}
\label{eq:linearform}
|n_i \log \delta-\log 2a-m_i\log \alpha|<\frac{9}{\alpha^{3m_i/2}}\qquad {\text{\rm holds for both}}\qquad  i=1,2,
\end{equation}
where $m_2>m_1>100$. We apply Matveev's theorem on the left--hand side of \eqref{eq:multform}.  First we need to check that
$$
\Gamma:=e^{\Lambda}-1= \delta^n (2a)^{-1} \alpha^{-m}-1
$$
is nonzero. Well, if it were, then $\delta^{n}=(2a) \alpha^{m}$. The right--hand side belongs to ${\mathbb Q}[\alpha]$ which is a field of degree $3$, while the left--hand side belongs to ${\mathbb Q}[{\sqrt{d}}]$ which is a quadratic field. The intersection of these two fields is 
${\mathbb Q}$. Hence, $\delta^{n}\in {\mathbb Q}$. Since $\delta$ is an algebraic integer and $n\ge 1$, it follows that $\delta^{n}\in {\mathbb Z}$. Since $\delta$ is a unit, we get that $\delta^{n}=1$, so $n=0$, a contradiction. Thus, $\Gamma\ne 0$, and we can apply Matveev's theorem. We take 
$$
l=3,~\eta_1=\delta,~\eta_2=2a,~\eta_3=\alpha,~d_1=n,~d_2=-1,~d_3=-m
$$
and ${\mathbb L}={\mathbb Q}[{\sqrt{d}},\alpha]$ which has degree $d_{\mathbb L}=6$. Since $\delta\ge 1+{\sqrt{2}}>\alpha$, the second inequality \eqref{eq:ttt} tells us right--away that $n<m$, so we take $D=m$.  We have $h(\eta_1)=(1/2)\log\delta$ and
$h(\eta_3)=(1/3)\log \alpha$. Further, 
$$
a=\frac{\alpha}{\alpha^2+2\alpha+3}
$$
and the minimal polynomial of $2a$ is $11X^3+4X-2$ and has roots $2a,~2b,~2c$. Further, $\max\{|2a|, |2b|, |2c|\}<1$ by \eqref{eq:num}. Thus, $h(\eta_2)=(1/3) \log 11$. Thus, we can take
$$
A_1=3\log\delta,\qquad A_2=2\log 11,\qquad A_3=2\log 1.84.
$$
Now Theorem \ref{thm:Matveev} tells us that
\begin{eqnarray*}
\label{eq:resultMatveev}
\log |\Gamma| & > & -1.4\times 30^6\times 3^{4.5}\times 6^2 (1+\log 6)(1+\log m) (3\log \delta) (2\log 11)(2\log 1.84 \nonumber\\
& > & -2.6\times 10^{14} \log \delta (1+\log m).
\end{eqnarray*}
Comparing the above inequality with \eqref{eq:multform}, we get
$$
1.5 m\log \alpha-\log 4.5<2.6\times 10^{14} \log \delta (1+\log m).
$$
Thus,
$$
m\log \alpha<1.8 \times 10^{14} \log \delta(1+\log m).
$$
Since $\alpha^m>\delta^{n}$ (see the second equation \eqref{eq:ttt}), we get that
\begin{equation}
\label{eq:neversusm}
n<1.8 \times 10^{14} (1+\log m).
\end{equation}
Further, since $\alpha>1.83$, we get
\begin{equation}
\label{eq:mversusdelta}
m<3\times 10^{14} \log \delta (1+\log m).
\end{equation}
Let us record what we have proved so far.
\begin{lemma}
\label{lem:prel}
If $X_n=T_m$ and $m>100$, then 
$$
n<1.8\times 10^{14} (1+\log m)\qquad {\text{and}}\qquad m<3\times 10^{14} \log \delta (1+\log m).
$$
\end{lemma}
Next, we return to the two inequalities given by \eqref{eq:linearform}. Multiply the one for $i=1$ with $n_2$ and the one for $i=2$ with $n_1$, subtract them and apply the triangle inequality to get that
\begin{eqnarray}
\label{eq:2logs}
|(n_2-n_1)\log 2a+(n_2m_1-n_1m_2)\log\alpha| & = & |n_2(n_1\log \delta-\log 2a-m_1\log \alpha)\nonumber\\
& - & n_1(n_2\log \delta-\log 2a-m_2\log \alpha|\nonumber\\
& \le &  n_2|n_1\log\delta-\log 2a-m_1\log \alpha| \nonumber\\
& + & n_1|n_2\log \delta-\log 2a-m_2\log \alpha|\nonumber\\
& \le & \frac{9n_2}{\alpha^{3m_1/2}}+\frac{9n_1}{\alpha^{3m_2/2}}\nonumber\\
& < & \frac{18n_2}{\alpha^{3m_1/2}}.
\end{eqnarray}
We are all set to apply Theorem \ref{thm:LMN} with 
$$
l=2,\quad \eta_1=2a,\quad \eta_2=\alpha,\quad d_1=n_2-n_1,\quad d_2=n_2m_1-m_2n_1.
$$
The fact that $\eta_1$ and $\eta_2$ are multiplicatively independent follows because the norm of $\eta_1$ is $2/11$ while $\eta_2$ is a unit.  Observe that $n_2-n_1<n_2$, while by the absolute value inequality in \eqref{eq:2logs}, we have
$$
|n_2m_1-n_1m_2|\le (n_2-n_1) \frac{|\log 2a|}{\log \alpha}+\frac{12 n_2}{\alpha^{3m_1/2} \log \alpha}<2n_2,
$$
because $m_1>100$. We have ${\mathbb L}={\mathbb Q}[\alpha]$ which has $d_{\mathbb L}=3$. So, we can take
$$
\log B_1=\max\left\{h(\eta_1), \frac{|\log \eta_1|}{3}, \frac{1}{3}\right\}=\frac{\log 11}{3}
$$
and 
$$
\log B_2=\max\left\{h(\eta_2), \frac{\log \eta_2}{3},\frac{1}{3}\right\}=\frac{1}{3}.
$$
Thus,
$$
b'=\frac{(n_2-n_1)}{3\times (1/3)}+\frac{|n_2 m_1-n_1m_2|}{3\times (\log(11)/3)}<2n_2.
$$
Now Theorem \ref{thm:LMN} tells us that with
$$
\Lambda:=(n_2-n_1)\log 2a+(n_2m_1-n_1m_2)\log\alpha,
$$
we have
$$
\log |\Lambda|>-24.34\times 3^4 \max\{\log(2n_2)+0.14,7\}^2 \cdot (1/3) \cdot (\log (11)/3).
$$
Thus, 
$$
\log |\Lambda|>-526 \left(\max\left\{\log 2n_2+0.14,7\right\}\right)^2.
$$
Combining this with \eqref{eq:2logs}, we get
$$
1.5 m_1\log \alpha-\log(18n_2)<526 \left(\max\{\log (2n_2)+0.14, 7\}\right)^2.
$$
If $\log(2n_2)+0.14\le 7$, then $n_2\le 476$. The above inequality then gives
$$
1.5 m_1\log \alpha<526\times 7^2+\log(12\times 476),
$$
which gives $m_1\le 28444$. Hence, $n_1<n_2\le 476$ and $m_1\le 28444$ in this case. Assume next that $n_2>476$. Then
$$
1.5 m_1\log \alpha<526 (\log(2n_2)+0.14)^2+\log(18n_2)<528 (1+\log n_2)^2,
$$
which gives
\begin{equation}
\label{eq:m1}
m_1<583 (1+\log n_2)^2.
\end{equation}
Since $\alpha^{m_1}>\delta^{n_1}\ge \delta$ (see the second relation \eqref{eq:ttt}), we get
$$
\log \delta<m_1\log \alpha<356 (1+\log n_2)^2.
$$
Combining this with the second inequality of Lemma \ref{lem:prel} with $(n,m)=(n_2,m_2)$, together with the fact that $n_2<m_2$, we get
$$
m_2<3\times 10^{14} \times 356 (1+\log m_2)^3,
$$
giving $m_2<1.6\times 10^{22}$. Inserting this into the first inequality of Lemma \ref{lem:prel}, we get $n_2<10^{16}$, which together with \eqref{eq:m1} gives $m_1<835000$. Let us summarize what we have proved.

\begin{lemma}\label{le:jb0}
If $X_{n_i}=T_{m_i}$ for $i=1,2$ with $m_1<m_2$ (so $n_1<n_2$), then 
$$
m_1<835000,\qquad n_2<10^{16},\qquad m_2<1.6\times 10^{22}.
$$
\end{lemma}

To lower these bounds we use continued fractions on \eqref{eq:2logs}, and Baker-Davenport reduction on \eqref{eq:linearform}. 

\section{The final computations}\label{sec:6}

Put $\chi=-\log 2a/\log\alpha$. Inequality (\ref{eq:2logs}) implies

\begin{equation}
\label{eq:2logsnew}
\left|(n_2-n_1)\chi-(n_2m_1-n_1m_2)\right| <\frac{18n_2}{\alpha^{3m_1/2}\log\alpha}.
\end{equation}
Since 
\begin{equation}\label{eq:jb1}
\frac{18n_2}{\alpha^{3m_1/2}\log\alpha}<\frac{1}{2(n_2-n_1)},
\end{equation}
it follows that $(n_2m_1-n_1m_2)/(n_2-n_1)$ is a convergent of $-(\log 2a)/(\log\alpha)$. Indeed,
$\log\alpha<0.61$ and $m_1>100$, together with Lemma \ref{le:jb0} induce
\begin{equation}
\alpha^{3m_1/2}>6\cdot10^{33}>60n_2^2>60(n_2-n_1)n_2>\frac{36}{\log\alpha}(n_2-n_1)n_2,
\end{equation} 
which immediately leads to (\ref{eq:jb1}).

Obvoiusly, $n_2-n_1<n_2<10^{16}$. Let $[a_0,a_1,a_2,\dots]=[1,1,1,1,6,1,1,22,1,\dots]$ be the continued fraction expansion of $\chi$, and let $p_k/q_k$ be its $k^{th}$ convergent.
After a computer calculation we found that
$$
4999601640630812=q_{33}<10^{16}<24351826693265967=q_{34},
$$
further the maximum of $a_i$ $(i=0,1,\dots,34)$ is $22=a_7$. 
Hence,
$$
\frac{1}{24n_2}<\frac{1}{24(n_2-n_1)}<\left|(n_2-n_1)\chi-(n_2m_1-n_1m_2)\right| <\frac{18n_2}{\alpha^{3m_1/2}\log\alpha},
$$
and comparing the leftmost and rightmost expressions, by Lemma \ref{le:jb0} it gives $m_1\le 87.8$. Since we assumed that $m_1>100$, we conclude that $m_1\le 100$. Now \eqref{eq:mn} gives $n_1<69.2$.

These upper bounds (on $n_1$ and $m_1$) make it possible to compute all existing $n_1$ and $m_1$. Defining
\begin{eqnarray*}
P^+_n(X)&=&\frac{(X+{\sqrt{X^2-1}})^n+(X-{\sqrt{X^2-1}})^n}{2}\quad{\rm and}\quad \\
P^-_n(X)&=&\frac{(X+{\sqrt{X^2+1}})^n+(X-{\sqrt{X^2+1}})^n}{2},
\end{eqnarray*}
a computer search on the equations
$$
P^+_{n_1}(X_1)=T_{m_1}\quad{\rm and}\quad P^-_{n_1}(X_1)=T_{m_1}
$$
with $1\le m_1\le 100$ and $1\le n_1\le 69$, where $n_1<m_1$ results in only the following possibilities:

Besides the trivial case $n_1=1$ (for both equations), which implies $X_1=T_{m_1}$, the only nontrivial solutions are
$$
(n_1,m_1,X_1)=(2,5,2),\qquad{\rm and}\qquad (n_1,m_1,X_1)=(3,5,1)
$$
in the first, and in the second case, respectively.

The non-trivial solutions lead to $(d,Y_1)=(3,1)$, and $(d,Y_1)=(2,1)$, respectively.
Now, applying (\ref{eq:linearform}) and Lemma \ref{lem:Baker-Davenport} we determine all the solutions to equation \eqref{eq4}.
First observe, that
\begin{equation*}
\left|n_2\frac{\log\delta}{\log\alpha}-m_2+\chi\right|<\frac{9}{\alpha^{3/2m_2}\log\alpha}<14.8\cdot2.4^{-m_2}.
\end{equation*} 

Put $\delta_1=2+\sqrt{3}$ and $\delta_2=1+\sqrt{2}$. 
Taking the continued fraction expansion of $\log\delta_i/\log\alpha$ ($i=1,2$), such that the suitable denominator of it exceeds $6\cdot 10^{16}$, we found that
$$
q_{1,31}=156827205418169727\approx1.56\cdot10^{17},$$
and
$$q_{2,28}=98827474195551603\approx9.88\cdot10^{16}$$
is satisfactory for $i=1$ and $i=2$, respectively.
We now apply Lemma \ref{lem:Baker-Davenport}, with $m=n_2$, $n=m_2$, $k=m_2$, $A=14.8$, $B=2.4$, $M=10^{16}$, $\kappa=\log\delta_i/\log\alpha$ and $\mu=\chi$. Further, according to the two cases $Q=q_{1,31}$ and  $Q=q_{2,28}$, we get $\xi_1>0.039$ and $\xi_2>0.071$. Consequently, $m_2<49.9$, $n_2<23.1$ in the first case, and $m_2<48.7$, $n_2<33.7$ in the second case. However, since we assumed that $m_2>100$, we get a contradiction, so $m_2\le 100$ leading to 
$n_2\le 69.2$. Checking the last range we only obtained the possibilities:
$$
X_1=2=T_3\qquad{\rm and} \qquad X_2=7=T_5, 
$$   
and 
$$
X_1=1=T_1=T_2\qquad{\rm and}\qquad X_3=7=T_5, 
$$
respectively.

Finally, in order to check the trivial cases $n_1=1$, $X_1=T_{m_1}$, we used a brute force algorithm which essentially coincides the treatment of the non-trivial cases. For any $1\le m_1\le 100$ we determined the decomposition $T_{m_1}^2-\varepsilon=dY_1^2$, where $d$ is squarefree. In this way we find $\delta_{m_1}=X_1+\sqrt{d}Y_1$. Then we consider the first convergents of the continued fraction expansions of
\begin{equation}\label{jun15}
\frac{\log \delta_{m_1}}{\log \alpha},
\end{equation} 
such that the denominator is larger than $M=6\cdot 10^{16}$, and the $\xi$ value in Lemma \ref{lem:Baker-Davenport} is positive. The upper bounds on $m_2$ are always less than 100, which contradicts the assertion $m_2>100$. Thus only the cases $m_2\le100$ remain to verify. 
As conclusion, the trivial cases do not yield further solutions to (\ref{eq4}). 

To illustrate the treatment, take
$\varepsilon=1$, $m_1=17$. Now $T_{17}=10609$, $T_{17}^2-1=112550880=7034430\cdot4^2$, therefore $\delta_{m_{1}}=10609+4\sqrt{7034430}$.
The first denominator of the continued fractions corresponding to (\ref{jun15}), which is larger then $M$ is $q_{29}$, but the first denominator with positive $\xi$ is $q_{31}$ ($\xi>0.276$).
Lemma \ref{lem:Baker-Davenport} implies $m_2\le50$. However, since we assumed that $m_2>100$, we get $m_2\le 100$. But the equations $P_{n}^{\pm} (X)=T_m$ were already solved for $m\le 100$, so we get no further solutions.  

\section{Acknowledgements}

We thank the referee for comments which improved the quality of this paper. The work on this paper started when the last three authors visited School of Mathematics of the Wits University in May 2016. They thank this Institution for support.  The second author was partially supported by PAPIIT IN114016 and CONACyT project 219827. All authors  also thank Kruger Park for excellent working conditions.


\begin{thebibliography}{9999}

\bibitem{BMS} Y. Bugeaud, M. Maurice, and S. Siksek,  \textit{Classical and modular approaches to exponential Diophantine equations I. Fibonacci and Lucas perfect powers}, Annals of Mathematics, \textbf{163} (2006), 969--1018.

\bibitem{BL} J. J. Bravo and F. Luca,  \textit{On a conjecture about repdigits in k-generalized Fibonacci sequences}, Publ. Math. Debrecen \textbf{82} (2013),  623--639.

\bibitem{DLT} A. Dossavi-Yovo, F. Luca and A. Togb\'e, \textit{On the $x$-coordinate of Pell equations which are repdigits},
Publ. Math. Debrecen \textbf{88} (2016), 381--399.

\bibitem{dujella} A. Dujella, A. Peth\H o, \textit{A generalization of a theorem of Baker and Davenport}, Quart. J. Math. Oxford Ser. (2) \textbf{49} (1998), 291--306.

\bibitem{FL} B. Faye and F. Luca, \textit{On x-coordinates of Pell equations which are repdigits}, Preprint 2016. 

\bibitem{LMN} M. Laurent, M. Mignotte and Yu. Nesterenko, \textit{Formes lin\'eaires en deux logarithmes et d\'eterminants d'interpolation},  J. Number Theory \textbf{55} (1995), 285--321.

\bibitem{LT} F. Luca, A. Togb\' e, \textit{On the $x$-coordinates of Pell equations which are Fibonacci numbers}, Mathematica Scandinavica, to appear.

\bibitem{Matveev} E. M. Matveev, \textit{An explicit lower bound for a homogeneous rational linear form in logarithms of algebraic numbers, II}, Izv. Ross. Akad. Nauk Ser. Mat. \textbf{64} (2000), 125--180. English translation in Izv. Math. \textbf{64} (2000), 1217--1269.

\bibitem{S} W. R. Spickerman, \textit{Binet's formula for the Tribonacci numbers}, The Fibonacci Quarterly \textbf{20} (1982), 118--120.





\end{thebibliography}
\end{document}